\documentclass[10pt]{amsart}

\newlength{\basicwidth}
\newlength{\shortbasicwidth}
\newlength{\basicheight}
\numberwithin{equation}{section}

\begin{document}

\vspace{0.4cm}
\begin{center}
\title{Identities for combinatorial sums involving trigonometric functions}
\maketitle
\end{center}

\vspace{0.3cm}
\begin{center}
HORST ALZER$^a$ \quad\mbox{and} \quad SEMYON YAKUBOVICH$^b$
\end{center}

\vspace{0.7cm}
\begin{center}
$^a$ Morsbacher Stra\ss{e} 10, 51545 Waldbr\"ol, Germany\\
\emph{Email:} {\tt{h.alzer@gmx.de}}
\end{center}

\vspace{0.4cm}
\begin{center}
$^b$ Department of Mathematics, Faculty of Sciences, University of Porto, \\ Campo Alegre st. 687, 4169-007 Porto, Portugal \\
\emph{Email:} {\tt{syakubov@fc.up.pt}}
\end{center}

\vspace{2.cm}
{\bf{Abstract.}}  Let
$$
A_{m,n}(a)=\sum_{j=0}^m (-4)^j {m+j\choose 2j}\sum_{k=0}^{n-1} \sin(a+2k\pi/n) \cos^{2j}(a+2k\pi/n)
$$
and
$$
B_{m,n}(a)=\sum_{j=0}^m (-4)^j {m+j+1\choose 2j+1}\sum_{k=0}^{n-1} \sin(a+2k\pi/n) \cos^{2j+1}(a+2k\pi/n),
$$
where $m\geq 0$ and $n\geq 1$ are integers and $a$ is a real number. We present two proofs for the following results:\\
(i) \, {If} $2m+1 \equiv 0 \,  (\mbox{mod} \, n)$, {then}
$$
A_{m,n}(a)=(-1)^m n \sin((2m+1)a).
$$
(ii) \, {If} $2m+1 \not\equiv 0 \, (\mbox{mod} \, n)$, {then} $A_{m,n}(a)=0$.

(iii)  \, {If} $2(m+1) \equiv 0 \,  (\mbox{mod} \, n)$, {then}
$$
B_{m,n}(a)=(-1)^m \frac{n}{2} \sin(2(m+1)a).
$$
(iv) \, {If} $2(m+1) \not\equiv 0 \, (\mbox{mod} \, n)$, {then} $B_{m,n}(a)=0$.


\vspace{1.3cm}

\vspace{0.01cm}

\newpage

\vspace{2cm}

\section{Introduction and statement of the main  result}

The work on this note has been inspired by a remarkable paper on trigonometric identities published by
 Chu and Marini \cite{CM} in 1999. One of their results states that
\begin{equation}
\sum_{k=0}^{n-1}\frac{\sin(a+2k\pi/n)}{1+x^2-2x\cos(a+2k\pi/n)}=G_{n,a}(x)
\end{equation}
with
\begin{equation}
G_{n,a}(x)=\frac{nx^{n-1}\sin(na)}{1+x^{2n}-2x^n \cos(na)}.
\end{equation}
Some authors used (1.1) and several similar identities to derive numerous interesting trigonometric summation formulas. We refer to
 Andrica and Piticari \cite{AP}, Chamberland \cite{Cha}, Chen \cite{Che}, Chu \cite{Chu}, Wang and Zheng \cite{WZ},
and the references cited therein.

Here, we study the  two closely related combinatorial sums
$$
A_{m,n}(a)=\sum_{j=0}^m (-4)^j {m+j\choose 2j}\sum_{k=0}^{n-1} \sin(a+2k\pi/n) \cos^{2j}(a+2k\pi/n)
$$
and
$$
B_{m,n}(a)=\sum_{j=0}^m (-4)^j {m+j+1\choose 2j+1}\sum_{k=0}^{n-1} \sin(a+2k\pi/n) \cos^{2j+1}(a+2k\pi/n),
$$
where $m\geq 0$ and $n\geq 1$ are integers and $a$ is a real number. We show that (1.1) and (1.2) can be applied to find identities for
$A_{m,n}(a)$ and $B_{m,n}(a)$. Our main result reads as follows.

\vspace{0.4cm}
{\bf{Theorem.}} \emph{Let $m\geq 0$, $n\geq 1$ be integers and let $a$ be a real number.}\\
(i) \emph{If} $2m+1 \equiv 0 \,  (\mbox{mod} \, n)$, \emph{then}
$$
A_{m,n}(a)=(-1)^m n \sin((2m+1)a).
$$
(ii) \emph{If} $2m+1 \not\equiv 0 \, (\mbox{mod} \, n)$, \emph{then} $A_{m,n}(a)=0$.

(iii) \emph{If} $2(m+1) \equiv 0 \,  (\mbox{mod} \, n)$, \emph{then}
$$
B_{m,n}(a)=(-1)^m \frac{n}{2} \sin(2(m+1)a).
$$
(iv) \emph{If} $2(m+1) \not\equiv 0 \, (\mbox{mod} \, n)$, \emph{then} $B_{m,n}(a)=0$.

\vspace{0.4cm}
These formulas can be used to deduce further identities. For example, if we integrate $A_{m,n}(a)$ from $a=0$ to $a=\pi/2$, then we obtain
\begin{displaymath}
C_{m,n}= \left\{ \begin{array}{ll}
(-1)^m {n}/(2m+1), & \textrm{if \,  $2m+1 \equiv 0 \,  (\mbox{mod} \, n)$, }\\
0, & \textrm{if \,  $2m+1 \not\equiv 0 \,  (\mbox{mod} \, n)$, }
\end{array} \right.
\end{displaymath}
where
$$
C_{m,n}= \sum_{j=0}^m \frac{(-4)^j}{2j+1} {m+j\choose 2j} \sum_{k=0}^{n-1} \bigl( \sin^{2j+1}(2k\pi/n)+\cos^{2j+1}(2k\pi/n)\bigl).
$$

In the following sections, we present two proofs of our Theorem. The first one makes use of a comparison of coefficients in series expansions, where (1.1) and (1.2)  play a key role. The second proof applies properties of the classical Chebyshev polynomials.

\vspace{0.3cm}
\section{First proof}

Let
$$
z=x(x-2\cos(a+2k\pi/n)).
$$
Then, for $|x|<\sqrt{2}-1$,
$$
|z|\leq |x| ( |x|+2)<1.
$$
Using (1.1), the geometric series and the binomial formula yields
\begin{eqnarray}
G_{n,a}(x) & = & \sum_{k=0}^{n-1}\sin(a+2k\pi/n)\sum_{\nu=0}^\infty (-x)^{\nu} (x-2\cos(a+2k\pi/n))^{\nu} \\ \nonumber
& = &  \sum_{\nu=0}^\infty \sum_{k=0}^{n-1}\sin(a+2k\pi/n) (-x)^{\nu}\sum_{j=0}^{\nu}
{\nu\choose j} x^j (-2\cos(a+2k\pi/n))^{\nu-j} \\ \nonumber
& = &  \sum_{\nu=0}^\infty \sum_{j=0}^{\nu} S(j,\nu) \nonumber
\end{eqnarray}
with
$$
S(j,\nu)=(-1)^j {\nu\choose j}   2^{\nu-j}  x^{\nu+j}\sum_{k=0}^{n-1}\sin(a+2k\pi/n) \cos^{\nu-j}(a+2k\pi/n).
$$
Then,
\begin{eqnarray}
\sum_{j=0}^m S(j,2m-j)
& = & x^{2m} \sum_{j=0}^m (-1)^j {2m-j\choose j} 2^{2(m-j)}\sum_{k=0}^{n-1} \sin(a+2k\pi/n) \cos^{2(m-j)}(a+2k\pi/n) \\ \nonumber
& = & (-1)^m A_{m,n}(a) x^{2m} \\ \nonumber
\end{eqnarray}
and
\begin{eqnarray}
\sum_{j=0}^m S(j,2m+1-j)
& = & x^{2m+1}\sum_{j=0}^m (-1)^j {2m+1-j\choose j} 2^{2(m-j)+1} \\ \nonumber
& & \times \sum_{k=0}^{n-1} \sin(a+2k\pi/n) \cos^{2(m-j)+1}(a+2k\pi/n) \\ \nonumber
& = &  2 (-1)^m B_{m,n}(a)  x^{2m+1}. \nonumber
\end{eqnarray}
Applying
$$
\sum_{\nu=0}^\infty \sum_{j=0}^\nu S(j,\nu)=\sum_{m=0}^\infty \sum_{j=0}^m S(j,2m-j)+\sum_{m=0}^\infty \sum_{j=0}^m
S(j,2m+1-j)
$$
we conclude from (2.1), (2.2) and (2.3) that
$$
G_{n,a}(x)=\sum_{m=0}^\infty \frac{G^{(2m)}_{n,a}(0)}{(2m)!}x^{2m}
+\sum_{m=0}^\infty \frac{G^{(2m+1)}_{n,a}(0)}{(2m+1)!}x^{2m+1}=
\sum_{m=0}^\infty  (-1)^m A_{m,n}(a)        x^{2m}+\sum_{m=0}^\infty   2 (-1)^m B_{m,n}(a)      x^{2m+1}.
$$
Comparing the coefficients reveals that
\begin{equation}
A_{m,n}(a)=(-1)^m \frac{G^{(2m)}_{n,a}(0)}{(2m)!}
\end{equation}
and
\begin{equation}
B_{m,n}(a)=(-1)^m \frac{G^{(2m+1)}_{n,a}(0)}{2\cdot (2m+1)!}.
\end{equation}

Next, we use the formula
$$
\sum_{k=1}^\infty \sin(k\alpha) z^k = \frac{z \sin(\alpha)}{1+z^2-2z\cos(\alpha)};
$$
see Prudnikov et al. \cite[Entry 5.4.12.1]{PBM}.
We set $z=x^n$, $\alpha=na$, multiply both sides by $n/x$ and use (1.2).  This gives
$$
\sum_{k=1}^\infty n\sin(kna) x^{nk-1}= G_{n,a}(x).
$$
It follows that for nonnegative integers $p$ we obtain
\begin{equation}
G^{(p)}_{n,a}(0)=0, \quad\mbox{if} \quad p+1\not\in n\mathbb{Z},
\end{equation}
and
\begin{equation}
G^{(p)}_{n,a}(0)=p! n \sin(kna),  \quad\mbox{if} \quad p=nk-1.
\end{equation}
We consider four cases.

{\underline{Case 1.}} \, $2m+1 \equiv 0 \, (\mbox{mod} \, n)$. \\
Let $n k_0=2m+1$. Using (2.4) and (2.7) gives
$$
A_{m,n}(a)=(-1)^m n \sin(k_0 na)=(-1)^m n \sin((2m+1)a).
$$

{\underline{Case 2.}} \, $2m+1 \not\equiv 0 \, (\mbox{mod} \, n)$. \\
From (2.4) and (2.6) we conclude that $A_{m,n}(a)=0$.

{\underline{Case 3.}} \, $2(m+1) \equiv 0 \, (\mbox{mod} \, n)$. \\
Let $n k_1=(2m+1)+1$. We apply (2.5) and (2.7). This yields
$$
B_{m,n}(a)=(-1)^m\frac{n}{2} \sin(n k_1 a)
=(-1)^m\frac{n}{2} \sin(2(m+1)a).
$$

{\underline{Case 4.}} \, $2(m+1) \not\equiv 0 \, (\mbox{mod} \, n)$. \\
Then, $(2m+1)+1\not\in n \mathbb{Z}$, so that (2.5) and (2.6) lead to $B_{m,n}(a)=0$.

The proof of the Theorem is complete.

\vspace{0.3cm}
\section{Second proof}

The  Chebyshev polynomials of the first and second kind  are defined by
$$
T_n(x)=\cos(n\theta) \quad\mbox{and} \quad U_n(x)=\frac{1}{n+1} T'_{n+1}(x)= \frac{\sin((n+1)\theta)}{\sin(\theta)},
\quad n=0,1,2,...,
$$
where $\cos(\theta)=x$. Detailed information about these functions can be found, for example, in Mason and Handscomb \cite{MH}.

In what follows, we make use of the following known formulas:
\begin{equation}
T_{2n+1}(\sin (\theta))=(-1)^n \sin((2n+1)\theta),
\end{equation}
\begin{equation}
\cos(\theta) U_{2n+1}(\sin(\theta))=(-1)^n \sin(2(n+1)\theta).
\end{equation}
Moreover, we need the trigonometric identities
\begin{equation}
\cos(  (2n+1)\arcsin(\cos(\theta))) =(-1)^n  \text{sign}(\sin(\theta)) \sin((2n+1)\theta),
\end{equation}
\begin{equation}
\sin(2n\arcsin(\cos(\theta)))=(-1)^{n+1}\text{sign}(\sin(\theta)) \sin(2n\theta)
\end{equation}
and the summation formula
\begin{equation}
\sum_{k=0}^N \sin(y+kz)= \frac{\sin(y+Nz/2)   \sin((N+1)z/2)}{\sin(z/2)}.
\end{equation}

\vspace{0.2cm}
We are now in a position to present our second proof of the Theorem. First, we show that (i) and (ii) are valid.
Let $x\in (-2,2)$. The following formula is given in Prudnikov et al. \cite[Entry 4.2.3.15]{PBM}:
\begin{equation}
\sum_{j=0}^m (-1)^j {m+j\choose 2j} x^{2j}=\frac{\cos((2m+1)\arcsin(x/2))}{\cos(\arcsin(x/2))}.
\end{equation}
We apply  (3.3) and (3.1). Then,
\begin{equation}
\frac{\cos(  (2m+1)\arcsin(  \cos(\theta))    )}{\cos(\arcsin(\cos(\theta)))}=\frac{ T_{2m+1}(\sin(\theta))}{\sin(\theta)}.
\end{equation}
Using  (3.6) with $x=2\cos(  a+2k\pi/n)$, (3.7)  and (3.1) leads to
\begin{equation}
A_{m,n}(a)  =   \sum_{k=0}^{n-1} T_{2m+1}(\sin(a+2k\pi/n))
 =  (-1)^m \sum_{k=0}^{n-1} \sin ( (2m+1) (a+2k\pi/n) ).
\end{equation}
If  $2m+1 \equiv 0 \,  (\mbox{mod} \, n)$, then we conclude from (3.8) that
$$
A_{m,n}(a)=(-1)^m n \sin((2m+1)a).
$$
Let  $2m+1 \not\equiv 0 \,  (\mbox{mod} \, n)$. We apply  (3.8) and (3.5) with $N=n-1$, $y=(2m+1)a$, $z=(2m+1)2\pi/n$. Then we obtain
$$
A_{m,n}(a)= (-1)^{m+1} \frac{\sin((2m+1)\pi)}{\sin((2m+1)\pi/n)}\sin( (2m+1)(a-\pi/n))=0.
$$

 Next, we prove (iii) and (iv). We have
\begin{equation}
\sum_{j=0}^m (-1)^j {m+j+1\choose 2j+1} x^{2j+1}=\frac{\sin(2(m+1)\arcsin(x/2))}{\cos(\arcsin(x/2))};
\end{equation}
see Prudnikov et al. \cite[Entry 4.2.3.17]{PBM}. Using (3.4), (3.3) and (3.2)  gives
\begin{equation}
\frac{\sin(  2(m+1)\arcsin(\cos(\theta))    )}{\cos(\arcsin(\cos(\theta)))}=\frac{\cos(\theta)  U_{2m+1}(\sin(\theta))}{\sin(\theta)}.
\end{equation}
We set $x=2\cos(a+2k\pi/n)$. Then we conclude from (3.9), (3.10) and (3.2) that
\begin{eqnarray}\nonumber
B_{m,n}(a) & = & \frac{1}{2} \sum_{k=0}^{n-1} \cos(a+2k\pi/n) U_{2m+1}(\sin(a+2k\pi/n)) \\ \nonumber
& = & \frac{(-1)^m}{2} \sum_{k=0}^{n-1} \sin ( 2(m+1) (a+2k\pi/n) ). \nonumber
\end{eqnarray}
If $2(m+1) \equiv 0 \,  (\mbox{mod} \, n)$, then we obtain
$$
B_{m,n}(a)=(-1)^m \frac{n}{2} \sin((2(m+1)a).
$$
Let $2(m+1) \not\equiv 0 \,  (\mbox{mod} \, n)$. Then, we conclude from (3.5) that
$$
B_{m,n}(a)=\frac{(-1)^m}{2} \frac{\sin(2(m+1)\pi)}{\sin( 2(m+1)\pi/n)}\sin( 2(m+1)(a-\pi/n))=0.
$$
This completes the proof of the Theorem.

\end{document}